\documentclass[oneside, a4paper,11pt,reqno]{amsart}
\usepackage{xcolor}
\usepackage{ulem}
\usepackage[utf8]{inputenc}

\newtheorem{hypo}{Hypothesis}

\newtheorem{thm}[hypo]{Theorem}

\newtheorem{lem}[hypo]{Lemma}

\def\PP{\mathbb{P}}
\def\RR{\mathbb{R}}
\def\ZZ{\mathbb{Z}}

\def\EE{\mathbb{E}}
\def\NN{\mathbb{N}}

\let\BFseries\bfseries\def\bfseries{\BFseries\mathversion{bold}} 

\def\deq{\stackrel{\stackrel{\mathcal{L}}{=}}{}}

\def \ind {{\bf 1}}

\def\dd{\mbox{d}}

\title{Persistence exponent for discrete-time, time-reversible processes}

 \author{ Frank Aurzada} 
\address{AG Stochastik, Fachbereich Mathematik, Technische Universit\"at Darmstadt, Schlossgartenstr.\ 7, 64289 Darmstadt, Germany.}
\email{aurzada@mathematik.tu-darmstadt.de}

\author{Nadine Guillotin-Plantard} 
\address{Institut Camille Jordan, CNRS UMR 5208, Universit\'e de Lyon, Universit\'e Lyon 1, 43, Boulevard du 11 novembre 1918, 69622 Villeurbanne, France.}
\email{nadine.guillotin@univ-lyon1.fr}

\subjclass[2010]{60G50, 60G22, 60G10, 60G15, 60F10, 62M10}
\keywords{First passage time, Fractional Gaussian noise, Long-range dependence, Persistence, Random walk, Random scenery. \\
This research was supported by the french ANR project MEMEMO2 2010 BLAN 0125}

\date{\today}

\begin{document}

\begin{abstract}
We study the persistence probability for some discrete-time, time-reversible processes.  In particular, we deduce the persistence exponent in a number of examples: first, we deal with random walks in random sceneries (RWRS) in any dimension with Gaussian scenery. Second, we deal with sums of stationary Gaussian sequences with correlations exhibiting long-range dependence. Apart from the persistence probability we deal with the position of the maximum and the time spent on the positive half-axis by the process.
\end{abstract}

\maketitle

\section{Introduction and main results}
\subsection{Introduction}
Persistence concerns the probability that a stochastic process has a long negative excursion. In this paper, we are concerned with discrete-time processes. If $Z=(Z_k)_{k=0,1,2,\ldots}$ is a stochastic process, we study the rate of the probability
$$
\PP\left[ \max_{k=1,\ldots, T} Z_k \leq 0 \right],\qquad \text{as $T\to+\infty$}.
$$
In many cases of interest, the above probability decreases polynomially, i.e., as $T^{-\theta+o(1)}$, and it is the first goal to find the persistence exponent $\theta$. For a recent overview on this subject, we refer to the survey \cite{AS} and for the relevance in theoretical physics we recommend \cite{Maj1,BMS13}.

The purpose of this paper is to analyse the persistence probability for time-reversible processes, i.e.\ processes such that, for any $T\geq 0$,
$$
(Z_{T-k}-Z_T)_{k=0,\ldots,T} \deq (Z_k)_{k=0,\ldots,T},
$$
where $\deq$ means equality in law. Note that this property implies in particular that $Z_0=0$ and that the $(Z_k)$ are symmetric. Another consequence of this property is that the increments of $Z$
 are stationary (see Lemma 2.1 in \cite{castellguillotinwatbled}).
 
The main message of this paper is that the time-reversibility property together with the ``natural scaling'' of the process already yield the persistence exponent. Of course, we will have to impose some other technical conditions to make our theorems work. One of them is Gaussianity, which helps a lot on a technical level (due to the use of Slepian's inequality, \cite{slepian}), but we do not think that this is the essential assumption. In particular, our first example will be a non-Gaussian process: random walk in random sceneries.

The main tool in the proofs is a relation of the persistence probability to the ``exponential functional''
\begin{equation} \label{eqn:expfunction}
\Phi(T):=\EE\left[ \Big(\displaystyle\sum_{l=0}^{[T]} e^{Z_l}\Big)^{-1}\right],\qquad T\geq 0,
\end{equation}
which is in some sense a smoothed out version of the persistence probability: for paths satisfying $\max_{k=1,\ldots,T} Z_k \leq 0$ the quantity $(\sum_{l=0}^{[T]} e^{Z_l})^{-1}$ will typically be reasonably large, while for $\max_{k=1,\ldots,T} Z_k > 0$ one can expect it to be relatively small. A formal relation between some exponential functional and the persistence probability was already discovered in the continuous-time setup, see e.g.\ \cite{Molchan1999}. This paper is the first rigorous treatment of the discrete-time counterpart. The quantity $\Phi$ above is also of own importance in theoretical physics, see e.g.\ \cite{oshaninetal}.

The paper is structured as follows. In Section~\ref{sec:resultsrwrs} we give the main results for random walks in random sceneries, while Section~\ref{sec:resultssumsof} contains the results for sums of stationary sequences with long range dependence, and in Section~\ref{sec:resultsfurther} we give a few generalizations. The result for exponential functionals is stated and proved in Section~\ref{sec:expfunctionals} and may be of independent interest. Sections~\ref{sec:rwrs}, \ref{sec:sumsof}, and~\ref{sec:prooffurther} contain the proofs for RWRS, sums of stationary sequences, and the generalizations, respectively.

\subsection{Random walks in random scenery} \label{sec:resultsrwrs}
Random walks in random sceneries were introduced independently by H. Kesten and F. Spitzer \cite{KS} and by A. N. Borodin \cite{Borodin}.
Let $S = (S_n)_{n \ge 0}$ be a random walk in $\mathbb{Z}^d$ starting at $0$,
i.e., $S_0 = 0$ and
$
X_n:= S_n-S_{n-1}, n \ge 1, \text{ is a sequence of i.i.d.\ } \mathbb{Z}^d \text{-valued random variables}.
$
Let $\xi = (\xi_x)_{x \in \mathbb{Z}^d}$ be a field of i.i.d.\ real random variables independent of $S$.
The field $\xi$ is called the random scenery.
{\it The random walk in random scenery (RWRS)} $Z := (Z_n)_{n \ge 0}$ is defined 
by setting $Z_0 := 0$ and, for $n \in \mathbb{N}^{*}$,
\begin{equation}
Z_n := \sum_{i=1}^n \xi_{S_i}.
\end{equation}
We will denote by $\mathbb{P}$ the joint law of $S$ and $\xi$. Limit theorems for RWRS have a long history, we refer to \cite{GuPo} for a complete review.

In the following, we consider the case when the random scenery $(\xi_x)_{x\in \ZZ^d}$  is assumed to be Gaussian with zero mean and variance equal to one. We are interested in the persistence properties
of the sum $Z_n, n\geq 1$. 

We distinguish two cases: When $d=1$ and the walk is recurrent and in the domain of attraction of an $\alpha$-stable law, the scaling limit of RWRS is a continuous-time non-Gaussian process, called Kesten-Spitzer process, \cite{KS}. Contrary, when $d\geq 2$ or the random walk is transient, the scaling limit of RWRS is the Brownian motion \cite{Bo89,S76,KS}.

{\bf Case $d=1$:} Let us first consider the one-dimensional case. Concerning the random walk $(S_n)_{n\geq1}$, the distribution of $X_1$ 
is assumed to be centered and to belong to the normal domain of attraction of a strictly stable distribution  $\mathcal{S}_{\alpha}$ of 
index $\alpha \in (1,2]$, with characteristic function $\phi$ given by
$$\phi(u)=e^{-|u|^\alpha(A_1+iA_2 \text{sgn}(u))}\quad u\in\mathbb{R},$$
where $0<A_1<\infty$ and $|A_1^{-1}A_2|\le |\tan (\pi\alpha/2)|$. Under the previous assumptions, the following weak convergence holds in the space  of 
c\`adl\`ag real-valued functions 
defined on $[0,\infty)$,  endowed with the 
Skorokhod topology:
$$\left(n^{-\frac{1}{\alpha}} S_{\lfloor nt\rfloor}\right)_{t\geq 0}   
\mathop{\Longrightarrow}_{n\rightarrow\infty}
^{\mathcal{L}} \left(Y(t)\right)_{t\geq 0},$$
where $Y$ is a L\'evy process such 
that $Y(0)=0$, 
$Y(1)$ has distribution $\mathcal{S}_{\alpha}$.
We will denote by $(L_t(x))_{x\in\mathbb{R},t\geq 0}$ a continuous version with compact support of the local time of the process
 $(Y(t))_{t\geq 0}$, \cite{marcusrosen}.
In \cite{KS}, Kesten and Spitzer proved 
the convergence in distribution of $((n^{-(1-1/2\alpha)} Z_{[nt]})_{t\ge 0})_n$, 
to a process $\Delta=(\Delta_t)_{t\geq 0}$ defined by
$$\Delta_t : = \int_{\mathbb{R}} L_t(x) \, \dd W(x),$$
where $(W(x))_{x\ge 0}$ and $(W(-x))_{x\ge 0}$ are independent standard Brownian motions independent of $Y$. The process $\Delta$ is called Kesten-Spitzer process in the literature.
Our main result in this setup is the following.
\begin{thm}
\label{theoMk}
Consider RWRS $(Z_k)$ for Gaussian scenery with $d=1$, the random walk being in the normal domain of attraction of an $\alpha$-stable law (as detailed above). Then there exists a constant $c>0$ such that for large enough $T$
\begin{equation}
\label{eqMk}
 T^{-1/2\alpha} e^{- c \sqrt{\log T} }\leq \PP\Big[ \max_{k=1,\ldots,T} Z_k \leq 0\Big]
\leq  T^{-1/2\alpha} (\log T)^{+c}.
\end{equation}
\end{thm}

The corresponding results for the continuous-time Kesten-Spitzer process $\Delta$ were obtained in \cite{BFFN}, also cf.\ \cite{Maj,MdM,castellguillotinwatbled}. We stress that the proofs of the present results do not follow from \cite{BFFN}, however.

{\bf Case $d=2$ and transient random walks:}
When $d=2$ we assume that the random walk increment $X_1$ has a centered law with a finite and non-singular covariance matrix $\Sigma$.
We further suppose that the random walk is aperiodic in the sense of Spitzer \cite{S76},
which amounts to requiring that $\varphi(u) = 1$ if and only if $u \in 2 \pi \ZZ^2$, where $\varphi$ is the characteristic function of $S_1$.
In this situation, Bolthausen \cite{Bo89} proved that a functional central limit theorem is satisfied for $((n\log n)^{-1/2}Z_{[nt]})_{t\geq0})$,  namely
$$\left( (n\log n)^{-1/2} Z_{[ nt ]}\right)_{t\geq 0}   \mathop{\Longrightarrow}_{n\rightarrow\infty}^{\mathcal{L}} \left(\sigma B(t)\right)_{t\geq 0},$$
where $(B(t))_{t\geq 0}$ is a one dimensional Brownian motion and  $\sigma^2 = (\pi \sqrt{\det \Sigma})^{-1}$.\\*
The case of transient random walks has been dealt with in Spitzer \cite{S76} and Kesten and Spitzer \cite{KS}: Here,
$$\left( n^{-1/2} Z_{[ nt ]}\right)_{t\geq 0}   \mathop{\Longrightarrow}_{n\rightarrow\infty}^{\mathcal{L}} \left(\sigma B(t)\right)_{t\geq 0},$$
where $\sigma^2=2 G(0,0) - 1$ and $G$ is the Green function of the random walk. Note that this includes random walks evolving in $\ZZ^d, d\geq 3$, as well as
any transient random walk in $\ZZ^d, d=1,2$.

\begin{thm}\label{theoMkd}
Consider RWRS $(Z_k)$ for Gaussian scenery either with $d=2$ and the above assumptions on the random walk or with a transient random walk. Then there exists a constant $c>0$  such that for large enough $T$
\begin{equation} \label{eqn:thmd2}
T^{-1/2}  e^{- c \sqrt{\log T} }  \leq \PP\Big[ \max_{k=1,..,T} Z_k \leq 0\Big]
\leq  T^{-1/2} (\log T)^{+c}.
\end{equation}
\end{thm}

The proof of Theorem~\ref{theoMk} is given in Section~\ref{sec:rwrs}. The proof of Theorem \ref{theoMkd} is omitted since it is very similar to the one of Theorem \ref{theoMk}.

\subsection{Sums of stationary sequences}  \label{sec:resultssumsof}
Let $(X_i)_{i\geq 0}$ be a stationary Gaussian sequence with mean $0$ and variance $1$, and with correlations 
$r(j):=\EE [X_0 X_j] = \EE[ X_k X_{j+k} ]\geq 0$ satisfying as $n\rightarrow +\infty$,
\begin{equation} \label{eqn:lrd}
\sum_{i,j=1}^n r(i-j) \sim K  n^{2H} \ell(n), 
\end{equation}
where $H\in [\frac{1}{2},1)$, $K>0$, and $\ell$ is a slowly varying function at infinity. We are interested in the persistence exponent of the sum of the stationary sequence $Z_n:=\sum_{i=1}^n X_i$ for $n\geq 1$ and $Z_0:=0$. Recall that the scaling limit of $(Z_n)$ is the fractional Brownian motion with Hurst parameter $H$, \cite{taqqu}:
$$
\left( n^{-H} \ell(n)^{-1/2} Z_{[ nt ]}\right)_{t\geq 0}   \mathop{\Longrightarrow}_{n\rightarrow\infty}^{\mathcal{L}} \left( B^H(t)\right)_{t\geq 0}.
$$
A stationary sequence satisfying relation (\ref{eqn:lrd}) is said to have {\it long-range dependence}. We refer to \cite{samorodnitsky} for a recent overview of the field.

In this setup we can show the following theorem.

\begin{thm} \label{thm:statseq}
Under the above assumptions on the sequence $(X_i)$, in particular (\ref{eqn:lrd}) with $H\geq 1/2$, there is some constant $c>0$ such that for large enough $T$
\begin{equation} \label{eqn:thmlrd}
T^{-(1-H)}  \ell(T) e^{- c \sqrt{\log T} } \leq \PP\left[ \max_{k=1,\ldots,T} Z_k \leq 0 \right] \leq T^{-(1-H)} \ell(T) (\log T)^{+c}.
\end{equation}
\end{thm}

A simple example is fractional Gaussian noise: If $B^H$ denotes a fractional Brownian motion, define $X_i:=B^H_{i+1}-B^H_{i}$, $i=0,1,\ldots$. If $H\geq 1/2$ then Theorem~\ref{thm:statseq} holds and gives
$$
T^{-(1-H)} (\log T)^{-c} \leq \PP\left[ \max_{k=1,\ldots,T} B^H_k \leq 0 \right] \leq T^{-(1-H)} (\log T)^{+c}.
$$
The better lower bound follows from \cite{Aurzada}, while the upper bound is new.

We stress that we can prove Theorem~\ref{thm:statseq} for $H\ge 1/2$ only. We conjecture that it also holds for $H<1/2$.

We remark that the value $0$ of the boundary in the persistence probability is of no importance in the present results. In fact, in the situation of Theorems~\ref{theoMk}, \ref{theoMkd}, and \ref{thm:statseq}, we have
\begin{equation} \label{eqn:1stepslepian}
\PP\left[ \max_{k=1,\ldots,T} Z_k \leq a \right] \geq \PP[ \mathcal{N}(0,1)\le a-b] \cdot \PP\left[ \max_{k=1,\ldots,T-1} Z_k \leq b \right],
\end{equation}
for any $a\in\mathbb{R}, b\geq 0$; so that for any $a\in\mathbb{R}$ the terms $\PP\left[ \max_{k=1,\ldots,T} Z_k \leq a \right]$ and $\PP\left[ \max_{k=1,\ldots,T} Z_k \leq 0 \right]$ differ at most by a multiplicative constant.

The proof of the last remark and of Theorem~\ref{thm:statseq} can be found in Section~\ref{sec:sumsof}.

\subsection{Quantities related to the persistence probability}  \label{sec:resultsfurther}
Let us finally consider two terms that are closely related to the persistence probability. Define the position of the maximum of the first $T$ steps of the process $Z$:
$$
\tau_T := \arg\max \{ Z_k, k=0,\ldots, T\}
$$
and the time spent by $Z$ in the positive half-axis up to time $T$:
$$
N_T:=| \{ k \leq T : Z_k>0\} |.
$$

Recall that the persistence probability concerns a scenario where $Z$ has a long negative excursion. One would expect that this scenario happens typically whenever $\tau_T$ and $N_T$ are small. This can be concretized by the following fact.

\begin{thm} \label{thm:tauandN}
Let $Z$ be one of the processes from Theorems~\ref{theoMk}, \ref{theoMkd}, or \ref{thm:statseq}. Assume $n=n(T)$ is an $\mathbb{N}_+$-valued function with $n=T^{o(1)}$. Then
\begin{eqnarray*}
\PP\left[ \max_{k=1,\ldots,T} Z_k \leq 0 \right]~\leq &  \PP[ \tau_T < n] & \le ~ \PP\left[ \max_{k=1,\ldots,T} Z_k \leq 0 \right] \cdot T^{o(1)},
\\
\PP\left[ \max_{k=1,\ldots,T} Z_k \leq 0 \right] ~\le &\PP[ N_T < n]  & \le~ \PP\left[ \max_{k=1,\ldots,T} Z_k \leq 0 \right] \cdot T^{o(1)}.
\end{eqnarray*}
\end{thm}

\section{Relation to exponential functionals} \label{sec:expfunctionals}
The main idea is to relate the persistence probability to the exponential functional (\ref{eqn:expfunction}), cf.\ \cite{Molchan1999,Aurzada,BFFN,AurzadaBaumgartenjpa,castellguillotinwatbled}. In \cite{Molchan1999} it is shown that the continuous-time analog of this quantity behaves as $c T^{H-1}$ for any {\it continuous-time} $H$-self-similar process with stationary increments and a certain other time-reversibility property. Further, certain moment conditions are assumed in \cite{Molchan1999} (also see \cite{Molchan1999preprint}). Note that we do not require these moment conditions in this paper for the below Lemma~\ref{Mol} to hold. Further, $H$-self-similarity, which of course does not make sense in discrete time, is replaced by (\ref{ass1}), which extracts the ``natural scaling'' of the process $Z$.

The result is as follows.

\begin{lem}\label{Mol} Let $Z=(Z_n)_{n\in\NN}$ be a stochastic process with 
\begin{equation} \label{ass1}
\lim_{T\to+\infty} \frac{1}{T^{H}\ell(T)} \EE\left[ \sup_{  t\in[0,1]} Z_{[t T]} \right] = \kappa,
\end{equation}
for some $H\in(0,1)$, $\kappa\in(0,\infty)$, and with $\ell$ being a slowly varying function at infinity. Further assume that $Z$ is time-reversible in the sense that for any $T\in\NN$, the vectors $(Z_{T-k}-Z_T)_{k=0,\ldots,T}$ and $(Z_k)_{k=0,\ldots,T}$ have the same law. 
Then, 
$$\limsup_{x\rightarrow +\infty} \frac{x^{1-H}}{ \ell(x)} \EE\Big[ \Big(\displaystyle\sum_{l=0}^{[x]} e^{Z_l}\Big)^{-1} \Big]\leq \kappa H$$
and
$$\liminf_{x\rightarrow +\infty} \frac{x^{1-H}}{\ell(x)} \EE\Big[ \Big(\displaystyle\sum_{l=1}^{[x]} e^{Z_l}\Big)^{-1} \Big]\geq \kappa H.$$
\end{lem}

Note the difference in the summation $l=0,\ldots$ vs.\ $l=1,\ldots$, which complicates the use of this lemma. In particular, if one can show the second assertion with summation in $l=0,1,\ldots$, one would be able to reduce the error term in the lower bound from $e^{-c\sqrt{\log T}}$ to $(\log T)^{-c}$ in (\ref{eqMk}), (\ref{eqn:thmd2}), and (\ref{eqn:thmlrd}); cf.\ (\ref{eqn:summation}).

\begin{proof}
Let us define for every $T\in [1,+\infty),$
$$ \Psi(T):= \EE\left[ \log\Big(\sum_{k=0}^{[T]-1} e^{Z_k} + (T- [T])e^{Z_{[T]} } \Big)\right].$$
We clearly have
$$\EE\Big[ \sup_{t\in [0,1]} Z_{[t([T]-1)]}\Big] \leq \Psi(T) \leq \EE\Big[ \sup_{t\in [0,1]} Z_{[tT]}\Big]  +\log (T+ 1).$$
From assumption (\ref{ass1}), it follows that $\Psi(T)\sim \kappa T^H \ell(T) $.

By Fubini's theorem we have for any $u\in (1,+\infty)$,
\begin{eqnarray*}
\Psi(u) &=& \EE\Big[ \log \Big( \sum_{l=0}^{[u]-1} e^{Z_l} + (u-[u]) e^{Z_{[u]} } \Big) \Big] \\
&=& \int_1^u \Psi'(x) \dd x,
\end{eqnarray*}
where 
$$ \Psi'(x) =\sum_{k=1}^{\infty} \EE\left[ \frac{e^{Z_k}}{\sum_{l=0}^{k-1} e^{Z_l} +(x-k) e^{Z_k}} \right] {\bf 1}_{[k,k+1)}(x)=\sum_{k=1}^{\infty} \EE\left[ \frac{1}{\sum_{l=1}^{k} e^{Z_{k-l}-Z_k} +(x-k)} \right] {\bf 1}_{[k,k+1)}(x).$$
Using time reversibility,
$$ \Psi'(x) =\sum_{k=1}^{\infty} \EE\left[ \frac{1}{\sum_{l=1}^{k} e^{Z_l} +(x-k)} \right] {\bf 1}_{[k,k+1)}(x).$$
Let $0< a<b <+\infty$. Then,  for $x$ large enough,
\begin{eqnarray}\label{diff}
\Psi( [bx]+1) -\Psi( [ax]) &=& \int_{[ax]}^{[bx]+1} \Psi'(u) \dd u\nonumber \\
&=& \sum_{k=[ax]}^{[bx]}  \int_k^{k+1} \EE\left[ \frac{1}{\sum_{l=1}^{k} e^{Z_l} +(u-k)} \right] \, \dd u.
\end{eqnarray}
\noindent{\bf Estimation of the limsup:} Estimating the last quantity from below, we get the inequality:
\begin{eqnarray*}
\Psi ( [bx]+1) -\Psi( [ax])&\geq & ([bx]  -[ax] +1) \EE\Big[ \Big(\sum_{l=0}^{[bx]} e^{Z_l} \Big)^{-1}\Big].
\end{eqnarray*}
Therefore, 
$$ \frac{x^{1-H}}{\ell(x)} \EE\Big[ \Big(\sum_{l=0}^{[bx]} e^{Z_l}\Big)^{-1} \Big]\leq \frac{\Psi( [bx]+1) -\Psi( [ax])}{ x^H \ell(x)}\frac{x}{([bx]  -[ax] +1 )}.  $$
Since
\begin{eqnarray*}
\frac{\Psi( [bx]+1) -\Psi( [ax])}{ x^H \ell(x)} &=& \frac{\Psi( [bx]+1)}{ x^H \ell(x)} - \frac{\Psi( [ax])}{ x^H \ell(x)} \\
&=& \frac{\Psi( [bx]+1)}{ ([bx]+1)^H \ell([bx]+1) }\frac{([bx]+1)^H \ell([bx]+1)}{x^H \ell(x)}
\\
&&- \frac{\Psi( [ax])}{ [ax]^H \ell([ax])}\frac{[ax]^H \ell([ax])}{x^H \ell(x)}\\
&\rightarrow & \kappa (b^H -a^H),
\end{eqnarray*} 
as $x\rightarrow +\infty$, we obtain
$$\limsup_{x\rightarrow +\infty} \frac{x^{1-H}}{\ell(x)} \EE\Big[ (\displaystyle\sum_{l=0}^{[bx]} e^{Z_l})^{-1} \Big]\leq \frac{\kappa (b^H -a^H) }{(b-a)}.$$
Now taking $b=1$ and $a\uparrow 1$, we get 
$$\limsup_{x\rightarrow +\infty} \frac{x^{1-H}}{\ell(x)} \EE\Big[ (\displaystyle\sum_{l=0}^{[x]} e^{Z_l})^{-1} \Big]\leq \kappa H.$$
{\bf Estimation of the liminf:}  First, from (\ref{diff}),  we have the inequality
\begin{eqnarray*}
\Psi( [bx]+1) -\Psi( [ax])&\leq & ([bx]  -[ax] +1 ) \EE\Big[ \Big(\sum_{l=1}^{[ax]} e^{Z_l} \Big)^{-1}\Big]
\end{eqnarray*}
Therefore, 
$$ \frac{x^{1-H}}{\ell(x)}  \EE\Big[ \Big(\sum_{l=1}^{[ax]} e^{Z_l}\Big)^{-1} \Big]\geq \frac{\Psi( [bx]+1) -\Psi( [ax])}{ x^H \ell(x)}\frac{x}{([bx]  -[ax] +1)}  $$
Since
\begin{eqnarray*}
\frac{\Psi( [bx]+1) -\Psi( [ax])}{ x^H \ell(x)} &=& \frac{\Psi( [bx]+1)}{ x^H \ell(x)} - \frac{\Psi( [ax])}{ x^H \ell(x)} \\
&=& \frac{\Psi( [bx]+1)}{ ([bx]+1)^H  \ell([bx]+1)}\frac{([bx]+1)^H \ell([bx]+1)}{x^H \ell(x)}
\\
&& - \frac{\Psi( [ax])}{ [ax]^H \ell([ax])}\frac{[ax]^H \ell([ax])}{x^H \ell(x)}\\
&\rightarrow & \kappa (b^H -a^H),
\end{eqnarray*}
as $x\rightarrow +\infty$, we obtain
$$\liminf_{x\rightarrow +\infty} \frac{x^{1-H}}{\ell(x)} \EE\Big[ (\displaystyle\sum_{l=1}^{[ax]} e^{Z_l})^{-1} \Big]\geq \frac{\kappa (b^H -a^H) }{(b-a)}.$$
Now taking $a=1$ and $b\downarrow 1$, we get 
$$\liminf_{x\rightarrow +\infty} \frac{x^{1-H}}{\ell(x)}  \EE\Big[ (\displaystyle\sum_{l=1}^{[x]} e^{Z_l})^{-1} \Big]\geq \kappa H.$$
\end{proof}

\section{Proof of Theorem \ref{theoMk}} \label{sec:rwrs}
\subsection{Verification of Lemma~\ref{Mol} for RWRS}
The goal of this subsection is to verify that Lemma~\ref{Mol} holds with $H:=1-\frac 1{2\alpha}$ and $\ell\equiv 1$ when $Z$ is a RWRS. First, by conditioning on the random walk it is seen easily that $Z$ is time-reversible in the sense that for any $T\in\NN$ the vectors $(Z_{T-k}-Z_T)_{k=0,\ldots,T}$ and $(Z_k)_{k=0,\ldots,T}$ have the same law.

In order to verify (\ref{ass1}) let us define the self-intersection local time of the random walk $S$ by
$$V_n:= \sum_{i,j=1}^n {\bf 1}_{\{S_i =S_j\}}.$$

Note that the sequence of random variables $T^{-H}\max_{k=1,\ldots,T} Z_k$ is uniformly bounded in $L^2$: Indeed, given $S$, 
the random variable $Z_n$ is a sum of associated random variables with zero mean and finite variance, so from Theorem~2 in \cite{NW}, $$\EE[ \max_{k=1,\ldots, T} Z_k^2 | S] \leq \EE [ Z_{T}^2| S] = V_{T}.$$
By integrating with respect to the random walk, we get $$\EE[ \max_{k=1,\ldots, T} Z_k^2] \leq \EE [ V_{T}] \sim C T^{2H},$$
cf.\ (2.13) in \cite{KS}.
Since the sequence of processes  $( Z_{[tT]} / T^{H} )_{t\geq 0}$ weakly converges to the process $(\Delta_t)_{t\geq 0}$ (see \cite{KS}), we get 
$$ \lim_{T\rightarrow +\infty}\EE\Big[ \sup_{t\in [0,1]} \Big(\frac{Z_{[tT]}}{T^{H}}\Big)\Big] =
\EE\Big[ \sup_{t\in [0,1]} \Delta_t\Big] =:\kappa,$$
which is known to be finite using Proposition 2.1 in \cite{BFFN}.


\subsection{Proof of the upper bound} 

As in \cite{Molchan1999} and \cite{Aurzada}, the main idea in the proof of the upper
bound in \eqref{eqMk}, is to bound the exponential functionals treating in Lemma~\ref{Mol} (cf.\ (\ref{eqn:expfunction})) from below 
by restricting the expectation  to a well-chosen set of paths.

Let us denote by $N_n(x), x\in \ZZ$, the local time of the random walk $S$ up to time $n$ and let us rewrite
$$Z_n = \sum_{x\in \ZZ} N_n(x) \xi_x.$$
\noindent Conditionally on $S$, the process $(Z_k)_{k \ge 0}$
is a centered Gaussian process such that for every $0\le l<k$, 
$$\EE[Z_l Z_k|  S] = \sum_{x\in \ZZ} N_l(x) N_k(x)   \geq 0 , $$
$$ \EE[Z_l (Z_k - Z_l)|  S] = \sum_{x\in\ZZ} N_l(x) ( N_k(x) - N_l(x))  \geq 0 ,$$
since $k \mapsto N_k(x)$ is increasing for all $x \in {\mathbb Z}$. 
It follows then from Slepian's lemma \cite{slepian},  that for every $0\le u<v<w$ and all real numbers $a,b$, 
\begin{equation}\label{slepian1}
\PP\left[\max_{k=u,\ldots,v} Z_k\le a,\ \max_{k=v+1,\ldots,w}Z_k\le b\Big|  S \right]\ge  
    \PP\left[\max_{k=u,\ldots,v}Z_k\le a\Big|  S \right]
    \PP\left[\max_{k=v+1,\ldots,w}Z_k\le b\Big|  S \right]
\end{equation}
\begin{equation}\label{slepian2}
\PP\left[\max_{k=u,\ldots,v}Z_k\le a,\ \max_{k=v+1,\ldots,w}(Z_k-Z_v)\le b \Big|  S \right]\ge  
    \PP\left[\max_{k=u,\ldots,v}Z_k\le a \Big|  S \right]\PP\left[\max_{k=v+1,\ldots,w}(Z_k-Z_v)\le b\Big|  S\right].
\end{equation}

Let $a_T:= [\log T]^2$ and 
set $\beta_T:=\sqrt{V_{a_T}}$ where $V_{a_T}=\sum_{x\in\ZZ} N_{a_T}(x)^2$.
Let us define the random function 
$$\phi(k):= \left\{\begin{array}{ll}
  1  & \mbox{ for } 0 \leq k < a_T \, , \\
1-\beta_T  &  \mbox{ for } a_T\leq k \leq T \, ,
\end{array} 
\right.
$$ 
which is $S$-measurable.
Clearly, we have 
\begin{equation} \label{eqn:argumentphi1}
\EE\left[\left( \sum_{k=0}^T e^{Z_k}\right)^{-1}\Big| S \right] \geq 
\left( \sum_{k=0}^T e^{\phi(k)} \right)^{-1}  \PP\Big[\forall k\in\{0,\ldots,T\}, Z_k\leq \phi(k)\Big|  S \Big].
\end{equation}
By Slepian's lemma (see (\ref{slepian1})), we have
\[
 \PP\Big[\forall k\in\{0,\ldots,T\}, Z_k\leq \phi(k)\Big|  S \Big]
 \geq   \PP\Big[\max_{k=0,\ldots,a_T} Z_k\leq 1\Big|  S  \Big] \,  \PP\Big[\max_{ k=a_T,\ldots,T} 
Z_k\leq 1-\beta_T\Big|  S \Big].
\]
Note that
\begin{eqnarray*}
\PP\Big[\max_{k=a_T,\ldots,T} Z_k\leq 1-\beta_T \Big|  S \Big] 
& \geq & \PP\Big[Z_{a_T} 
\leq -\beta_T ; \max_{k=a_T,\ldots,T} (Z_k- Z_{a_T}) \leq 1\Big|  S \Big]
\\
& \geq & \PP\Big[Z_{a_T} \leq -\beta_T\Big| S \Big] \cdot \PP\Big[\max_{k=a_T,\ldots,T} 
(Z_k - Z_{a_T} )\leq 1\Big|  S \Big],
\end{eqnarray*}
by Slepian's lemma (see (\ref{slepian2})). Let $\overline{\Phi}$ be the function defined for every $ u\in \RR$ by
\begin{equation}\label{eqn:defnbarphi}
\overline{\Phi}(u)=\frac{1}{\sqrt{2\pi}} \int_u^{+\infty} e^{-\frac{x^2}{2} }\dd x,
\end{equation}
so that we have 
$$ \PP\left[Z_{a_T} \leq -\beta_T| S\right] =\overline{\Phi}(\beta_T V_{a_T}^{-1/2})= \overline{\Phi}(1).$$
Moreover, it is easy to check that for every $T>1$
$$\sum_{k=0}^T e^{\phi(k)}  \leq e  (a_T+1+T e^{-\beta_T}).$$
In the following, $C$ is a constant whose value may change but does not depend on $T$. 
Then, summing up (\ref{eqn:argumentphi1}) and the succeeding estimates, we can write that for $T$ large enough
\begin{equation}\label{firststep}
\EE\left[\left(\sum_{k=0}^T e^{Z_k} \right)^{-1}\Big| S \right]
\geq C(a_T+Te^ {-\beta_T})^{-1} \, 
\PP\big[\max_{k=0,\ldots,a_T}Z_k\leq 1\big| S \big]\,
\PP\big[\max_{k=a_T,\ldots,T}(Z_k-Z_{a_T})\leq 1 \big| S \big].
\end{equation}
Next we use the maximal inequality in Proposition 2.2 in \cite{KL}
to write
$$\PP\big[\max_{k=0,\ldots,a_T}Z_k\leq 1\big| S\big]
=1-\PP\big[\max_{k=0,\ldots,a_T}Z_k > 1\big| S \big]
\geq 1-2\PP\big[Z_{a_T}\geq 1\big| S \big]
=\PP\left[|Z|\leq V_{a_T}^{-1/2}| S\right],$$
where $Z$ is a Gaussian variable $\mathcal{N}(0,1)$ independent of $S$, from which we deduce that there
exists a constant $c>0$ such that 
\begin{equation}\label{morceau01}
\PP\big[\max_{k=0,\ldots,a_T}Z_k\leq 1\big| S \big]
\geq c \min( V_{a_T}^{-1/2} ,1).
\end{equation} 
Injecting \eqref{morceau01} into \eqref{firststep}
we get that for $T$ large enough,
$$\PP\left[\left.\max_{k=a_T,\ldots,T}(Z_k-Z_{a_T})\leq 1 \right|S \right]
\leq
C\EE\left[\left(\sum_{k=0}^T e^{Z_k} \right)^{-1}\Big| S \right]
(a_T+Te^ {-V_{a_T}^{1/2}})\,\max( V_{a_T}^{1/2},1) .$$
Since for every $n\geq 1$, $n\leq V_n\leq n^2$, we get by taking expectations of the last inequality and using Lemma~\ref{Mol}:
$$\begin{aligned}
\PP\big[\max_{k=a_T,\ldots,T} (Z_k-Z_{a_T})\leq 1\big] &\leq  C (\log T)^4\,  T^{-\frac{1}{2\alpha}}.
\end{aligned}$$
The left hand side is greater than the quantity we want to bound from above, since by stationarity of increments,
\begin{equation} \label{eqn:o}
\PP\big[\max_{k=a_T,\ldots,T}(Z_k-Z_{a_T})\leq 1\big]
=\PP\big[\max_{k=0,\ldots,T-a_T}Z_ {k}\leq 1\big]
\geq \PP\big[\max_{k=0,\ldots,T}Z_ {k}\leq 1\big].
\end{equation}

\subsection{Proof of the lower bound}

Fix $\beta>1-H$ and define $Z_T^*:=\max_{k=1,\ldots,T} Z_k$. Observe that
\begin{eqnarray} \label{eqn:argumentlowerboundalpha}
 \EE \Big[ \Big(\sum_{k=1}^T e^{Z_k} \Big)^{-1} \Big] &= & \EE \Big[ \Big(\sum_{k=1}^T e^{Z_k} \Big)^{-1}\ind_{ Z_T^* \geq \beta \log T} \Big] + 
 \EE \Big[ \Big(\sum_{k=1}^T e^{Z_k} \Big)^{-1}\ind_{ Z_T^* < \beta \log T} \Big] \\
&= : & I_1(T) + I_2(T).
\end{eqnarray}
First, we clearly have 
\begin{eqnarray}
I_1(T) &\leq &\EE [ e^{-Z_T^*} \ind_{ Z_T^* \geq \beta \log T}]  \notag
\\
&\leq & T^{-\beta}. \notag
\end{eqnarray}
Secondly, let us define the event $A:=\{ Z_1 \geq - \sqrt{4\beta \log T} \}$. Then, 
\begin{equation}\label{eqn:summation}
I_2(T) \leq e^{ \sqrt{4\beta \log T}} \PP\left[ Z_T^* < \beta \log T \right] + \EE \Big[ \Big(\sum_{k=1}^T e^{Z_k} \Big)^{-1}\ind_{A^{c}} \Big]. 
\end{equation}
Using the Cauchy-Schwarz inequality and the symmetry of $Z_1$, the second term on the right hand side can be bounded as follows 
$$
\EE \Big[ \Big(\sum_{k=1}^T e^{Z_k} \Big)^{-1}\ind_{A^{c}} \Big]^2 \leq \PP\Big[ Z_1 \geq  \sqrt{4\beta \log T}\Big] \  \EE[e^{2 Z_1}].
$$
Since $Z_1 $ is $\mathcal{N}(0,1)$,  
$$\PP\Big[ Z_1 \geq  \sqrt{4\beta \log T}\Big] =\frac{1}{\sqrt{2\pi} } \int_{\sqrt{4\beta \log T}}^{+\infty} e^{-u^2/ 2} \, \dd u\leq   T^{-2\beta}.$$
Since $T^{-\beta}$ is of lower order, we have shown that for $T$ large,
\begin{equation}\label{ZN}
\PP\left[Z_T^* \leq \beta \log T\right] \geq c T^{H-1} e^{-\sqrt{4 \beta \log T}}.
\end{equation}

Let $a_T:=[\log T]^2$. Note that by Slepian's lemma (cf.\ (\ref{slepian1}),(\ref{slepian2})),
\begin{eqnarray} 
\PP\left[\left. \max_{k=1,\ldots,T} Z_k\leq 1 \right|S\right] &\geq& \PP\left[ \max_{k=1,\ldots,a_T} Z_k\leq 1; Z_{a_T}\leq - \beta \log T;\notag \right.\\ && \qquad\left.\left. \max_{k=a_T+1,\ldots, T} Z_k - Z_{a_T} \leq \beta \log T\right|S\right]
\notag \\
&\geq & \PP\left[ \left. \max_{k=1,\ldots,a_T} Z_n\leq 1\right|S\right] \cdot \PP\left[Z_{a_T}\leq - \beta \log T|S\right] \notag \\ && \qquad \cdot \PP\left[ \left.\max_{k=a_T+1,\ldots, T} Z_k - Z_{a_T} \leq \beta \log T\right|S\right]\label{eqn:logtermslepian}
\end{eqnarray}
Further, let $Z$ be a Gaussian random variable $\mathcal{N}(0,1)$ independent of $S$,
\begin{eqnarray}
\PP[Z_{a_T}\leq - \beta \log T|S] &=& \PP\Big[ Z \leq \frac{- \beta \log T}{V_{a_T}^{1/2}}\Big|S\Big] \notag \\
&\geq & \overline{\Phi}\Big( \beta \frac{\log T}{[\log T]}\Big) \geq \overline{\Phi}\Big(2\beta\Big), \label{eqn:gaussiandirect}
\end{eqnarray}
since $V_{a_T}\geq a_T = [\log T]^2$ and where $\overline{\Phi}$ is defined in (\ref{eqn:defnbarphi}).

Thus, rewriting (\ref{eqn:logtermslepian}) and using the last inequality, inequality (\ref{morceau01}), and stationary increments as in (\ref{eqn:o}), we get 
\begin{eqnarray} \label{ZN2}
\PP\left[  \max_{k=1,\ldots,T}  Z_k \leq \beta \log T\right]
&\leq & \EE \left[\PP\left[ \max_{k=a_T+1,\ldots, T} Z_k - Z_{a_T} \leq \beta \log T|S\right]\right] \nonumber\\
&\leq & c \EE\left[\PP\left[ \left.\max_{k=1,\ldots,T} Z_k\leq 1\right|S\right]\cdot \max( V_{a_T}^{1/2},1) \right]\nonumber\\
&\leq & c (\log T)^2\,  \PP\left[ \max_{k=1,\ldots,T} Z_n\leq 1\right] 
\end{eqnarray}
having used that $V_{a_T}\leq a_T^2\leq (\log T)^4$.
By combining (\ref{ZN}) and (\ref{ZN2}), we get the lower bound.

\section{Proof of Theorem~\ref{thm:statseq}} \label{sec:sumsof}
\subsection{Verification of Lemma~\ref{Mol}.} \label{sec:proofofsumsof}
It is well-known (see for instance \cite{Pitt}) that positively correlated Gaussian random variables are associated, so from Theorem~2 in \cite{NW}, 
$$
\EE\Big[ \max_{k=1,\ldots, T} Z_k^2\Big] \leq \EE \Big[ Z_{T}^2\Big]  \sim K T^{2H}  \ell(T). 
$$
Thus, the sequence of random variables $(T^{2H} \ell(T))^{-1/2}\max_{k=1,\ldots, T} Z_k$ is uniformly bounded in $L^2$.
Moreover from Lemma 5.1 in \cite{taqqu} the sequence $(Z_{[nt]}/n^{H}\sqrt{\ell(n)})_{t\geq 0}$ weakly converges for the Skorokhod topology to 
$(\sqrt{K} B_H(t))_{t\geq 0}$ where $(B_H(t))_{t\geq 0}$ is the Fractional Brownian motion with Hurst index $H$. So assumption (\ref{ass1}) follows with 
$\kappa = \sqrt{K}\  \EE\big[ \displaystyle\sup_{t\in [0,1]} B_H(t)\big]$.

We remark that Lemma~\ref{Mol} also holds for sums of stationary sequences with $H<1/2$ (using Gaussian concentration rather than the present argument).

It is seen easily that $Z$ is time-reversible in law: the random vectors $(X_1,\ldots,X_T)$ and $(X_T,\ldots, X_1)$ have the same law (both are zero mean Gaussian random vectors
 with correlations $\EE[X_i X_j ] = \EE[ X_{T-(i-1)} X_{T-(j-1)}]$), so $(Z_0, Z_1,\ldots, Z_T)$ and $(Z_T- Z_{T-0},Z_{T} -Z_{T-1},\ldots, Z_T-Z_{T-T})$ have the same law, and by symmetry the required time-reversibility property follows.


\subsection{Proof of the upper bound.}
Define $a_T:=[\log T]^2$ and
$$
\phi(u):=\begin{cases} 1 & u \leq a_T,\\
		       -\beta \log T & u> a_T,
         \end{cases}
$$
where $\beta>1$. Then,
$$
\sum_{l=0}^T e^{\phi(k)} \leq c a_T.
$$
We obtain by the argument analogous to (\ref{eqn:argumentphi1}) that
$$
c' \ell(T) T^{H-1}  \geq (c a_T)^{-1} \PP\left[ \forall k\in\{1,\ldots,T\}, Z_k\leq \phi(k)\right].
$$
Using non-negative correlations and Slepian's lemma, we obtain
\begin{eqnarray*}
\PP\left[ \forall k\in\{1,\ldots,T\}, Z_k\leq \phi(k)\right] &=&  \PP\left[ \max_{ k=1,\ldots, a_T}Z_k \leq 1 ;  \max_{ k=a_T + 1,\ldots, T} Z_k \leq  -\beta \log T \right]
\\
&\geq& \PP\left[ \max_{ k=1,\ldots, a_T}Z_k \leq 1\right]\cdot \PP\left[ \max_{ k=a_T,\ldots, T} Z_k \leq  -\beta \log T\right].
\end{eqnarray*}
Observe that $Z_{a_T}$ is $\mathcal{N}(0,\sigma^2)$ with $\sigma^2 = \sum_{1\leq i,j \leq a_T} \EE[ X_i X_j]$. 
Using the fact that the correlations are non-negative and the Cauchy-Schwarz inequality, we get 
$$
a_T=\sum_{i=1}^{a_T} \EE[X_i^2] \leq  \sigma^2  \leq \sum_{i,j=1}^{a_T} \EE[X_0^2]= a_T^2.
$$
Applying the maximal inequality in Proposition 2.2 in \cite{KL}, we  can write
$$\PP\big[\max_{k=1,\ldots,a_T}Z_k\leq 1\big]
=1-\PP\big[\max_{k=1,\ldots,a_T}Z_k > 1\big]
\geq 1-2\PP\big[Z_{a_T} \geq 1\big]
=\PP[|Z|\leq \sigma^{-1}],$$
where $Z$ is a Gaussian variable $\mathcal{N}(0,1)$, from which we deduce that there
exists a constant $c>0$ such that 
\begin{equation}\label{morceau1}
\PP\big[\max_{k=1,\ldots,a_T}Z_k\leq 1\big]
\geq c \ a_T^{-1}.
\end{equation} 
Again, with non-negative correlations and Slepian's lemma we get
\begin{eqnarray*}
\PP\left[ \max_{ k=a_T,\ldots, T} Z_k \leq  -\beta \log T\right] &\geq& \PP\left[ Z_{a_T}\leq -\beta \log T; \max_{ k=a_T,\ldots, T} Z_k -Z_{a_T} \leq 0\right]
\\
&\geq &
\PP[ Z_{a_T}\leq -\beta \log T]\cdot\PP\left[\max_{ k=a_T,\ldots, T} Z_k -Z_{a_T} \leq 0 \right].
\end{eqnarray*}
Since $a_T=[\log T]^2$, we get
\begin{equation}
\PP[ Z_{a_T}\leq -\beta \log T] \geq \PP\Big[ Z \leq -\beta  \frac{\log T}{[\log T]}\Big] \geq  \PP\Big[ Z \geq 2 \beta \Big] = \mbox{const.}, \label{eqn:whystar}
\end{equation}
Putting the pieces together we have seen:
$$
\PP\left[\max_{k=1,\ldots,T} Z_k \leq 0\right]\leq \ell(T)\ T^{H-1} (\log T)^{c}.
$$

\subsection{Proof of the lower bound.}
The proof of the lower bound in the RWRS setup can easily be adapted: observe that $Z_1$ is $\mathcal{N}(0,1)$-distributed; note that one obtains (\ref{eqn:logtermslepian}) without any conditioning; analogously, (\ref{eqn:gaussiandirect}) holds with $V_{a_T}=\mathbb{V} Z_{a_T} = \sum_{i,j=1}^{a_T} \EE[ X_i X_j] \geq \sum_{i=1}^{a_T} \EE X_i^2 = a_T$; similarly, (\ref{ZN2}) transfers without the conditioning; and the derivation of (\ref{morceau01}) follows by the same argument and without conditioning.

\subsection{Proof of (\ref{eqn:1stepslepian}).}
Consider the case where $(Z_k)$ is a RWRS first. Then, using Slepian's lemma,
\begin{eqnarray*}
\PP\left[ \max_{k=1,\ldots,T} Z_k \leq a \right] &=& \EE\left[ \PP\left[ \left. \max_{k=1,\ldots,T} Z_k \leq a \right|S\right]\right]
\\
&\ge & \EE\left[ \PP\left[ \left. Z_1\leq a-b; \max_{k=2,\ldots,T} Z_k-Z_1 \leq b \right|S\right]\right]
\\
&\ge & \EE\left[ \PP\left[ \left. Z_1\leq a-b \right|S\right] \cdot \PP\left[ \left.\max_{k=2,\ldots,T} Z_k-Z_1 \leq b \right|S\right]\right].
\end{eqnarray*}
Using that $Z_1$ is standard normal, conditionally on $S$, the first term is a constant. Using the stationarity of increments, we get (\ref{eqn:1stepslepian}).

The case when $Z$ is a sum of a Gaussian stationary sequence is yet easier, as Slepian's lemma can be used directly.

\section{Proof of Theorem~\ref{thm:tauandN}} \label{sec:prooffurther}
\noindent {\bf Lower bound.} By the definition of $\tau_T$,
$$
\PP\left[ \max_{k=1,\ldots,T} Z_k \leq 0 \right]\leq \PP\left[ \max_{k=n,\ldots,T} Z_k \leq 0 \right] \leq \PP[ \tau_T < n]
$$
By the definition of $N_T$,
$$
\PP\left[ \max_{k=1,\ldots,T} Z_k \leq 0 \right]\leq \PP\left[ \max_{k=n,\ldots,T} Z_k \leq 0 \right] \leq \PP[ N_T < n].
$$
This already settles the lower bounds.

%

\noindent  {\bf Upper bound.} Consider the case where $(Z_k)$ is a RWRS. Fix $\gamma>1$ and $\beta>1-H$ (here $H$ is either $1-\frac{1}{2\alpha}$ or $\frac 1 2$ or $H$ according to the assumptions of Theorems~\ref{theoMk}, \ref{theoMkd}, or \ref{thm:statseq}) for the rest of the proof. First note that we may assume w.l.o.g.\ that $n\geq (2\beta\log T)^{1/2(\gamma-1)}$, as otherwise one may pass over to $\tilde n:=\max(n,[(2\beta\log T)^{1/2(\gamma-1)}]+1)$ and use that $\PP[ \tau_T < n ] \leq \PP[ \tau_T < \tilde n ]$.

For simplicity define again $Z_T^*:=\max_{k=1,\ldots,T} Z_k$. Fix $h:=n^\gamma$ with $\gamma>1$  and observe that on $\tau_T<n$ we have either $Z_T^*\leq h$ or $Z_n^*=Z_T^*>h$ and thus
$$
\PP[ \tau_T < n ] 
\leq \PP[ Z_T^* \leq h ] + \PP[ Z_{n}^* > h].
$$
Note that, from Proposition 2.2 in \cite{KL} and the fact that $V_{n}\leq n^2$,
\begin{eqnarray*}
\PP[ Z_{n}^* > h] &=& \EE\left[ \PP[ Z_{n}^* > h|S ]\right] \leq  2  \EE\left[ \PP[ Z_{n} > h|S ]\right]\\
&\leq &  \EE\left[ \exp\Big( - \frac{h^2}{2V_{n}}\Big)\right] \leq  \exp\Big( - {n}^{2(\gamma-1)}/2\Big)\leq  T^{-\beta},
\end{eqnarray*}
since $n\geq(2\beta\log T)^{1/2(\gamma-1)}$. On the other hand, we set $a_{n}:=h^2 = n^{2\gamma}$. By Slepian's lemma,
\begin{eqnarray}
& & \PP\left[\left. \max_{k=1,\ldots,T} Z_k \leq 1 \right|S\right] \notag
\\
&\ge& \PP\left[ \left. \max_{k=1,\ldots,a_{n}-1} Z_k \leq 1 ; Z_{a_{n}} \leq - h ; \max_{k=a_{n},\ldots,T} Z_k-Z_{a_{n}} \leq h\right|S \right] \notag
\\
&\ge& \PP\left[ \left. \max_{k=1,\ldots,a_{n}-1} Z_k \leq 1 \right|S\right]\cdot \PP[ Z_{a_{n}} \leq - h|S] \cdot \PP\left[ \left.\max_{k=a_{n},\ldots,T} Z_k-Z_{a_{n}} \leq h\right|S \right]. \label{eqn:--}
\end{eqnarray}
The first term is treated as in (\ref{morceau01}). The second term equals
$$
\PP[ Z_{a_{n}} \leq - h|S] = \PP[ \mathcal{N}(0,1) \leq - h V_{a_{n}}^{-1/2}|S] \geq \PP[ \mathcal{N}(0,1) \leq - h a_{n}^{-1/2}] = \text{const.}
$$
having used $V_{a_{n}}\geq a_{n}$. Using (\ref{eqn:--}), the last estimate, (\ref{morceau01}), and $V_{a_{n}}\leq a_{n}^2$, we obtain
\begin{eqnarray*}
\PP\left[ \max_{k=1,\ldots,T} Z_k \leq h \right]  &\leq & \PP\left[ \max_{k=a_{n},\ldots,T} Z_k - Z_{a_{n}} \leq h \right]
\\
&= &\EE\left[ \PP\left[ \left.\max_{k=a_{n},\ldots,T} Z_k - Z_{a_{n}} \leq h  \right|S \right]\right]
\\
 &\le &\EE\left[ \PP\left[ \left.\max_{k=1,\ldots,T} Z_k \leq 1\right|S \right] \cdot c \max(V_{a_{n}}^{1/2},1) \right]
\\
 &\le &\EE\left[ \PP\left[ \left.\max_{k=1,\ldots,T} Z_k \leq 1\right|S \right] \cdot c a_{n} \right].
\\
 &= & \PP\left[ \max_{k=1,\ldots,T} Z_k \leq 1 \right]\cdot c n^{2\gamma}.
\end{eqnarray*}

The proof of the upper bound for $\PP[ N_T < n]$ is similar and based on the following observations (cf.\ the proof of Theorem 2 in \cite{Molchan1999}):
$$
\PP[ N_T < n] \leq \PP[ Z_T^{*}\leq h] + \PP[ Z_T^{*}>h, N_T < n].
$$
We have just seen how to treat the first term as long as $h=T^{o(1)}$. The second term is handled with the following construction: Split the time interval $\{1,\ldots,T\}$ into intervals $(I_i)_{i=1,\ldots,M}$ of lengths between $n$ and $2n$ each. There are less than $[T/n]$ such intervals. If $N_T<n$, on every interval $I_i$ there has to be at least one $k_{0,i}$ with $Z_{k_{0,i}}\leq 0$ (as otherwise one would have $N_T\geq |I_i|\geq n$). If also $Z_T^*>h$ there must be one of these intervals with a fluctuation of $Z$ larger than $h$ (as $Z_T^*$ must be attained on one of the intervals). Using the stationarity of increments and the fact that the intervals are at most of length $2n$, one obtains
$$
\PP[ Z_T^*>h, N_T < n] \leq \sum_{i=1}^{M} \PP\left[ \max_{k\in I_i} Z_k - Z_{k_{0,i}} > h \right] \leq \sum_{i=1}^{[T/n]} \PP\left[ \max_{k=1,\ldots,2n} Z_k> h \right].
$$
The latter term can be treated as above when $n=T^{o(1)}$ and $h$ suitably chosen.

The case where $Z$ is a sum of stationary sequences can be treated in the same manner.


\begin{thebibliography}{00}
\bibitem{Aurzada} Aurzada, F. {\it On the one-sided exit problem for fractional Brownian motion.}
Electron. Commun. Probab., 16:392--404, 2011.

\bibitem{AurzadaBaumgartenjpa} Aurzada, F.; Baumgarten, C. {\it Persistence of fractional Brownian motion with moving boundaries and applications.}
Journal of Physics A: Mathematical and Theoretical 46 (2013), 125007.

\bibitem{AS} Aurzada, F.; Simon, T.  {\it Persistence probabilities \& exponents.} To appear in: L\'evy matters, Springer, arXiv:1203.6554, 2012.

\bibitem{Bo89}
Bolthausen, E. {\it A central limit theorem for two-dimensional random walks in random sceneries.}
Ann.\ Probab.\ 17 (1989) 108--115.



\bibitem{Borodin} Borodin, A. N. 
{\it A limit theorem for sums of independent random variables defined 
on a recurrent random walk.}
(Russian)  Dokl. Akad. Nauk SSSR 246(4):786--787, 1979. 

\bibitem{BMS13}
Bray, A.~J.; Majumdar, S.~N.; and Schehr, G.
\newblock {\it Persistence and first-passage properties in non-equilibrium systems.}
\newblock  Advances in Physics, 62(3):225--361, 2013.

\bibitem{BFFN} Castell, F.; Guillotin-Plantard, N.; P\`ene, F.; and Schapira, B. 
{\it On the one-sided exit problem for stable processes in random scenery.} 
Electron. Commun. Probab. 18(33):1--7, 2013.

\bibitem{castellguillotinwatbled} Castell, F.; Guillotin-Plantard, N.; and Watbled, F. {\it Persistence exponent for random processes in Brownian scenery}. Preprint, \verb+https://hal.archives-ouvertes.fr/hal-01017142v2+

\bibitem{Fe} Feller, W. An introduction to probability theory and its applications. Vol. II. Second edition, John Wiley and Sons, Inc., New York-London-Sydney, (1971). 

\bibitem{GuPo} Guillotin-Plantard, N. and Poisat, J. {\it  Quenched central limit theorems for random walks in random scenery}.
Stochastic Process. Appl. 123 (4) (2013) 1348--1367.


\bibitem{KS} Kesten, H. and Spitzer, F. 
{\it  A limit theorem related to a new class of self-similar processes.} 
Z. Wahrsch. Verw. Gebiete 50:5--25, 1979.  
%

\bibitem{KL} Khoshnevisan, D. and Lewis, T. M. {\it A law of iterated logarithm for stable processes in random scenery.}
Stochastic Process. Appl., 74(1):89--121, 1998.


\bibitem{Maj1} Majumdar, S. {\it Persistence in nonequilibrium systems}. Current Science 77 (3):370-375, 1999.

\bibitem{Maj}  Majumdar, S. {\it Persistence of a particle in the Matheron - de Marsily velocity field.} 
Phys. Rev. E 68, 050101(R), 2003.

\bibitem{marcusrosen}
Marcus, M. B. and Rosen, J. Markov processes, Gaussian processes, and local times. Cambridge Studies in Advanced Mathematics, 100. Cambridge University Press, Cambridge, 2006.

 \bibitem{MdM} Matheron, G. and de Marsily G.
 {\it Is transport in porous media always diffusive? A counterexample.}
 Water Resources Res. 16:901--907, 1980. 

\bibitem{Molchan1999preprint} Molchan, G.M. {\it Maximum of fractional Brownian motion: probabilities of small values.}
Preprint, \verb+https://www.ma.utexas.edu/mp_arc/c/00/00-195.ps.gz+


\bibitem{Molchan1999} Molchan, G.M. {\it Maximum of fractional Brownian motion: probabilities of small values.}
Comm. Math. Phys., 205(1):97--111, 1999.
%


%
%


\bibitem{NW} Newman, C. M. and Wright, A. L. {\it An invariance principle for certain dependent sequences}. Ann. Probab. 
9, (1981), no. 4, 671-- 675. 


\bibitem{oshaninetal}
Oshanin, G.; Rosso, A.; and Schehr, G. {\it Anomalous Fluctuations of Currents in Sinai-Type Random Chains with Strongly Correlated Disorder}.
Phys. Rev. Lett. 110 (2013), 100602.


\bibitem{Pitt} Pitt, L. {\it Positively correlated normal variables are associated}.  Ann.\ Probab.\ Vol. 10, No 2, (1982) 496 -- 499.

%
%

 \bibitem{samorodnitsky} Samorodnitsky, G. {\it Long range dependence.}
Found. Trends Stoch. Syst. 1 (2006), no. 3, 163--257.

\bibitem{slepian} Slepian, D.
{\it The one-sided barrier problem for Gaussian noise.}
Bell System Tech. J. 41 (1962), 463--501. 

\bibitem{S76}
Spitzer, F. Principles of Random Walks. Second ed., in: Graduate Texts in Mathematics, vol. 34, Springer-Verlag, New-York, 1976.

\bibitem{taqqu} Taqqu, M.S. {\it Weak convergence to fractional Brownian motion and to the Rosenblatt process.} Z. Wahrscheinlichkeitstheorie und Verw. Gebiete 31 (1974/75), 287--302.





%
%
%

\end{thebibliography}
\end{document}